\documentclass[12pt,leqno]{article}
\usepackage{amsmath,amsthm,amssymb}
\usepackage{verbatim}

\def\init{\setcounter{equation}{0}}
\newtheorem{theorem}{Theorem}[section]
\newcommand{\R}{\mathbb{R}}

\newtheorem{corollary}[theorem]{Corollary}

\usepackage{tikz}
\usetikzlibrary{decorations.pathreplacing}

\usepackage{verbatim}

\title{Remarks on the determination of the Lorentzian metric  by the lengths  of geodesics or null-geodesics}

\author{G.Eskin
 \ \ \  Department of Mathematics, UCLA,\\ Los Angeles,
CA 90095-1555, USA. \ E-mail: eskin@math.ucla.edu
}

\begin{document}

\maketitle

In blessed memory  of Misha Shubin

\begin{abstract}
We consider  a Lorentzian  metric  in $\R\times\R^n$.     
We show  that  if we know the lengths of the 
space-time  geodesics starting  at  $(0,y,\eta)$   when $t=0$,   then we can recover the metric at $y$.   
We prove the rigidity of Lorentzian metrics. 
We also  prove  a variant of  the rigidity  property 
for the case of null-geodesics: if two metrics  are close and if corresponding null-geodesics have equal Euclidian  lengths  
then the metrics are equal.
\end{abstract}

\section{Introduction}

Consider a Lorentzian  metric
\begin{equation}																				\label{eq:1.1}
\sum_{j,k=0}^n g_{jk}(x)dx_jdx_k
\end{equation}
in   $\R\times\R^n$,   where  $x_0\in \R $  is the time variable,  $x=(x_1,...,x_n)$  are  the space variables,   the metric  $[g_{jk}]_{j,k=0}^n$
is  
independent  of the  time variable  $x_0$  and the signature  of  the  matrix   $[g_{jk}]_{j,k=0}^n$  is  $(+,-,...,-)$.

Let
\begin{equation}																				\label{eq:1.2}
H(x,\xi_0,\xi)=\frac{1}{2}\sum_{j,k=0}^n g^{jk}(x)\xi_j\xi_k
\end{equation}
be  the corresponding Hamiltonian  where  $[g^{jk}(x)]_{j,k=0}^n$  is  the inverse  to  $[g_{jk}(x)]_{j,k=0}^n$.

Consider  the Hamiltonian  system
\begin{align}								        												     \label{eq:1.3}
&\frac{dx}{dt}=\frac{\partial H}{\partial \xi},\  \frac{d\xi}{dt}=-\frac{\partial H}{\partial x},\   x=(x_1,...,x_n),\ \xi=(\xi_1,...,\xi_n),			
\\
&\frac{dx_0}{dt}=\frac{\partial H}{\partial \xi_0}, \  \frac{d\xi_0}{dt}=-\frac{\partial H}{\partial x_0}=0,   			  \label{eq:1.4}
\\
\nonumber
\mbox{with initial conditions}
\\ & x_0(0)=0,\ x(0)=y,\ \xi(0)=\eta,\ \xi_0(0)=\eta_0.  															\label{eq:1.5}
\end{align}
We have
\begin{equation}																		\label{eq:1.6}
\frac{dH(x(t),\xi_0(t),\xi(t))}{dt}=\frac{\partial  H}{\partial x}\frac{\partial      H}{\partial \xi}+\frac{\partial H}{\partial \xi}\Big(-\frac{\partial H}{\partial x}\Big)
+\frac{\partial  H}{\partial\xi_0}\frac{d\xi_0}{dt}=0.
\end{equation}
Therefore
\begin{equation}																			\label{eq:1.7}
H(x(t),\xi_0(t),\xi(t))=H(y,\eta_0,\eta),\ \ \forall t.
\end{equation}
The solution  $x_0=x_0(t), x=x(t), \xi_0=\xi_0(t),\xi=\xi(t)$  of  (\ref{eq:1.3}),  (\ref{eq:1.4})  is called  a bicharacteristic  and its  projection  on  $(x_0,x)$
space-time is called geodesic.   When   $H(y,\eta_0,\eta)=0$ the curve $x_0=x_0(t), x=x(t),\xi_0=\xi_0(t),\xi=\xi(t)$ 
 is called  null-bicharacteristic  and it's
projection  on  $(x_0,x)$  space  is called  null space-time  geodesic.

Inverse  problems  of the recovery  of  the  Riemannian  metric  from
the lengths  
of  geodesics  were considered  in many  works  (see[1],  [3],  [4],  [5],  [6]  and others).
In  [7]  some  subclass of Lorentzian  metrics  was considered.

Now  we briefly  describe  the  content  of the paper.

The main  tool  of our approach  is the formula  (\ref{eq:2.5})  (see [1],  formula  (3.9)).   Using this formula  we prove
(Theorem \ref{theo:2.1})   that  if the  lengths of  geodesics  starting  at  any  point  $(0,y,\eta),\, y\in  \R^n, \eta\in \R^n$,  are  equal  then  the 
Lorentzian  metrics  equal  at  the point  $y\in\R^n$.

In \S4 and \S5  we prove  the rigidity  property  for  the Lorentzian  metrics.

\section{The determination of the metric}
\init
The following formula holds (see,  for example,  [1],  formula (3.9))
\begin{equation}												\label{eq:2.1}
(g(\xi_0,\xi),(\xi_0,\xi))=\Big(g^{-1}\Big(\frac{dx_0}{dt},\frac{dx}{dt}\Big),\Big(\frac{dx_0}{dt},\frac{dx}{dt}\Big)\Big),
\end{equation}
where  $g=[g^{jk}]_{j,k=0}^n$  is the same as  in  (\ref{eq:1.2}), 
$(g(\xi_0,\xi))_j=\sum_{k=0}^ng^{jk}\xi_k,\linebreak j=0,...,n,$
  $g^{-1}$  is the inverse of $g$.

To check  (\ref{eq:2.1})  we use  that  (see (\ref{eq:1.3}),  (\ref{eq:1.4}))
\begin{equation} 												\label{eq:2.2}
\Big(\frac{dx_0}{dt},\frac{dx}{dt}\Big)=g(\xi_0,\xi).
\end {equation}
Thus 
\begin{equation}	 											\label{eq:2.3}
(\xi_0,\xi)=g^{-1}\Big(\frac{d x_0}{dt},\frac{dx}{dt}\Big).
\end{equation}
Therefore
\begin{align}													\label{eq:2.4}
&(g(\xi_0,\xi),(\xi_0,\xi))=\Big(g\,g^{-1}\Big(\frac{dx_0}{dt},\frac{dx}{dt}\Big),g^{-1}\Big(\frac{dx_0}{dt},\frac{dx}{dt}\Big)\Big)
\\
\nonumber
=&\Big(\Big(\frac{dx_0}{dt},\frac{dx}{dt}\Big),g^{-1}\Big(\frac{dx_0}{dt},\frac{dx}{dt}\Big)
\end{align}
Thus  (\ref{eq:2.4})  is equivalent to (\ref{eq:2.1}).

Rewrite (\ref{eq:2.1})  using  (\ref{eq:1.7}).  We  get
\begin{align}												\label{eq:2.5}
&\sum_{j,k=0}^n  g_{jk}(x)\frac{dx_j}{dt}\frac{dx_k}{dt}=\sum_{j,k=0}^n g^{jk}(x)\xi_j\xi_k
\\
\nonumber
=& 2H(x(t),\xi_0(t),\xi(t))=2H(y,\eta_0,\eta).
\end{align}
Denote by $\Sigma^+$  the set where 
\begin{equation}												\label{eq:2.6}                                                                                                
H(y,\eta_0,\eta)>0.
\end{equation}

Therefore  on $\Sigma^+$  
we  have  
\begin{equation}												
\nonumber
\sum_{j,k=0}^n g_{jk}(x)\frac{dx_j}{dt}\frac{dx_k}{dt}>0.
\end{equation}

Let  $x_0=x_0(t,y,\hat\eta),x=x(t,y,\hat\eta)$  be 
the space-time  geodesic  starting  at  $(0,y,\hat\eta)$  when  $t=0$.  Here  $\hat\eta=(\eta_0,\eta)$.

Denoting  $\hat x=(x_0,x),\hat y=(0,y),$  we   can rewrite the  space-time  geodesic in a short form
$$
\hat x=\hat x(t,\hat y,\hat\eta).
$$

Let  $R(g,T,\hat y,\hat\eta)$   be the length  of the geodesic in $g^{-1}$  metric,  i.e.     
\begin{equation}												\label{eq:2.7}
R(g,T,\hat y,\hat\eta)=\int\limits_0^T \sqrt{\sum_{j,k=0}^n  g_{jk}(x)\frac{dx_j}{dt}\frac{dx_k}{dt}}dt.
\end{equation}
Using  (\ref{eq:2.5})  we get
\begin{equation}												\label{eq:2.8}
\sum_{j,k=0}^n  g_{jk}(x)\frac{dx_j}{dt}\frac{dx_k}{dt}=2H(y,\eta_0,\eta)=\sum_{j,k=0}^ng^{jk}(y)\eta_j\eta_k.
\end{equation}
Therefore
\begin{equation}												\label{eq:2.9}
R(g,T,\hat y,\hat\eta)=\int\limits_0^T \sqrt{\sum_{j,k=0}^n  g^{jk}(y)\eta_j\eta_k }\,dt
=\sqrt{2H(y,\eta_0,\eta)}T.
\end{equation}

\begin{theorem}												\label{theo:2.1}
Suppose  $g_1$  and  $g_2$  are two inverse  metric  tensors  in $\R\times\R^n$.
Suppose  integrals $ R(g_1,T,\hat y,\hat\eta)$  and  $R(g_2,T,\hat y,\hat\eta)$  are equal  for all  $\hat y,\hat\eta$.  Then  the metrics  $g_1^{-1}$   and  $g_2^{-1}$ 
 are also equal.
\end{theorem}

{\bf Proof.}  Since  $ R(g_1,T,\hat y,\hat\eta)=R(g_2,T,\hat y,\hat\eta)$   we have from   (\ref{eq:2.9})  that
$$
\Big( \sum_{j,k=0}^n  g_1^{jk}(y)\eta_j\eta_k \Big)^{\frac{1}{2} }T=
\Big( \sum_{j,k=0}^n  g_2^{jk}(y)\eta_j\eta_k \Big)^{\frac{1}{2}}T.
$$ 
Therefore
\begin{equation}												\label{eq:2.10}
g_1(y)=g_2(y).
\end{equation}
\qed

{\bf Remark  2.1}.   When  the  condition  (\ref{eq:2.6})  is  not satisfied,   i.e.  when
\begin{equation}												              \label{eq:2.11}                           
 \sum_{j,k=0}^n  g^{jk}(y)\eta_j\eta_k =0,
\end{equation}
the bicharacteristic  (\ref{eq:1.3}),  (\ref{eq:1.4})  is   the  null-bicharacteristic  and  its projection  on  $(x_0,x)$-space  is  the  
space-time null-geodesic.  

{\bf  Remark  2.2}.   
Consider the particular case        of the Riemannian metric  $[g_{jk}]_{j,k=1}^n$.
Thus  $g_{00}= g_{j0}=g_{0j}=0$.

Therefore formula (\ref{eq:2.9}) allows  to recover  $g'(y)$  as in Theorem \ref{theo:2.1}.

 When (\ref{eq:2.6})  holds,  $R(g,T,\hat y,\hat\eta)$  is the length of geodesics in $g^{-1}$  metric.
 When  (\ref{eq:2.11})  holds  the geodesic  is  the null-geodesic  and we are confronted  with the problem of defining
 the length of null-geodesics.   
   
Consider formula (2.8). 

If $H(y, \eta_0,\eta) =0$  then   
the length of  geodesic  in $g^{-1}$ metric is zero.
But the geodesic itself is not zero.   We called it  before a null-geodesic.
The definition of the length of geodesic must be modified  in the case of null-geodesic. 

 We   define  it  as
the Euclidian length $L$ of geodesic:    
\begin{align}
\nonumber														
L(g,T,y,\eta)=\int\limits_0^T\sqrt{\Big(\frac{dx_0}{dt}\Big)^2+|x_t|^2}dt,
\ \ \ |x_t|^2=\sum_{k=1}^n\Big|\frac{dx_k}{dt}\Big|^2,
\end{align}
where 
$x_0=x_0(t,y,\eta), x=x(t,y,\eta)$,   is the space-time   geodesics  
   starting at  $x_0=0,$     $(y,\eta)$  at  $t=0$  and
ending  when  $t=T$  at the point  $x_0(T),   x_T=x(T,y,\eta)$.

\section{Estimates  for the null-geodesics}
\init
Let  
$$
H_p=\frac{1}{2}\sum_{j,k=0}^n q_p^{jk}(x)\xi_j\xi_k=\frac{1}{2}q_p^{00}\xi_0^2
+\sum_{j=1}^n q_p^{0j}\xi_j\xi_0+\frac{1}{2}q_p' \xi\cdot \xi,
$$ 
$p=1,2,$ 
be two Hamiltonians.
Denote
\begin{equation}																\label{eq:3.1}
q=q_1+\tau(q_2-q_1),  \ 0\leq\tau\leq 1.
\end{equation}
Let  $x_\tau,\xi_\tau$  be  solution  of the  Hamiltonian 
system
\begin{align}																	\label{eq:3.2}
&\frac{dx_\tau}{dt}=q'(x_\tau(t))\xi_\tau(t),
\\
\nonumber
&\frac{d\xi_\tau}{dt}=-\frac{1}{2}\frac{\partial q'(x_\tau(t))}{\partial x}\xi_\tau(t)\cdot\xi_\tau(t),
\\
\nonumber
&x_\tau(t,y,\eta)\Big|_{t=0}=y,\ \ \xi_\tau(t,y,\eta)\Big|_{t=0}=\eta,
\end{align}
and  let  
$$
(3.2')\ \ \ \ \ \ \ \ \ \ \ \ \ \ \ \ \ \ \frac{dx_0^\tau}{dt}=q_\tau^{00}(x)\xi_0  +\sum_{j=1}^n  q_\tau^{0j} \xi_j,  \ \ x_0^\tau(0)=0. 
 \ \ \ \ \ \ \ \ \ \ \ \ \ \ \ \ \ \ \ \ \ \ \
$$
We shall study  the  behavior  of $(x_\tau(t,y,\eta),\xi_\tau(t,y,\eta)$  and  $x_0^\tau$  with respect to  $\tau$.
Differentiating  (\ref{eq:3.2})  in $\tau$  we get

\begin{equation}																\label{eq:3.3}
\frac{d}{dt}
\begin{pmatrix}
\frac{d x_\tau}{d\tau}
 \\  
\frac{d \xi_\tau}{d\tau}
\end{pmatrix}
=
Q
\begin{pmatrix}
\frac{d x_\tau}{d\tau}
 \\  
\frac{d \xi_\tau}{d\tau}
\end{pmatrix}
+
F
\end{equation}
where
\begin{align}																	\label{eq:3.4}
&Q=
\begin{bmatrix}
\frac{\partial q_1'}{\partial x}\xi_\tau &  q_1'\\
-\frac{1}{2}\big(\frac{\partial^2 q_1'}{\partial x^2}\xi_\tau\big)\cdot\xi_\tau
&- \frac{\partial  q_1'(x_\tau)\xi_\tau}{\partial  x}
\end{bmatrix}
, 
\\
\nonumber 
&F=
 \begin{bmatrix}
(q_2'-q_1')\xi_\tau +O(\tau(q_2'-q_1')^2)\xi_\tau
\\
-\frac{1}{2}\big(\big(\frac{\partial q_2'}{\partial x}-\frac{\partial q_1'}{\partial x}\big)\xi_\tau\big)\xi_\tau
+\big(O\big(\tau\big(\frac{\partial q_2'}{\partial x}-\frac{\partial q_1'}{\partial x}\big)^2\big)\xi_\tau\big)\cdot\xi_\tau
\end{bmatrix}.
\end{align}
Note that 
\begin{equation}														\label{eq:3.5}
\frac{dx_\tau}{d\tau}\Big|_{t=0}=0,\ \ \ \frac{d\xi_\tau}{d\tau}\Big|_{t=0}=0
\end{equation}
since  $x_\tau\big|_{t=0}=y,\ \ \xi_\tau\big|_{t=0}=\eta$.

We shall write   the solution  of the  Cauchy  problem  (\ref{eq:3.3}),  (\ref{eq:3.5})  in  the form
\begin{equation}														\label{eq:3.6}
\begin{bmatrix}
\frac{dx_\tau}{d\tau}
\\
\frac{d\xi_\tau}{d\tau}
\end{bmatrix}=R(t) F,
\end{equation}
where  $R(t)$  is the  solution  operator  of the  equation  (\ref{eq:3.3}).

If  $N$   is  large  enough  then  the following  
 estimate  
for  the solution  of the Cauchy  problem  (\ref{eq:3.3}),  (\ref{eq:3.5})  holds:
$$
\max_{0\leq t\leq T} e^{-Nt}  \Big(\Big|\frac{dx_\tau}{d\tau}\Big| +\Big|\frac{d\xi_\tau}{d\tau}\Big|              \Big)
\leq  C_N\int\limits_0^T e^{-Nt}  |F (x_\tau(t))|dt.
$$
Since  $q_2'-q_1'$  is bounded,   $\tau(q_2'-q_1')^2\leq  C|q_2'-q_1'|$. 
Thus
$
|F|\leq  C|q_2'-q_1'|+ C\big|\frac{\partial}{\partial  x}(q_2'-q_1')\big|.$
  Therefore
\begin{multline}														\label{eq:3.7}
\max_{0\leq t\leq T} e^{-Nt}  \Big(\Big|\frac{dx_\tau}{d\tau}\Big| +\Big|\frac{d\xi_\tau}{d\tau}\Big|              \Big)
\\
\leq C_N\,\sup_\tau\int\limits_0^T e^{-Nt}  \big|(q_2'-q_1') (x_\tau(t))\big|dt+C_N\,\sup_\tau\int\limits_0^T e^{-Nt}\Big|\frac{\partial(q_2'-q_1')}{\partial x}
(x_\tau(t))\Big|dt
\end{multline}
where here and below  $C_N$   means  various  constants  depending  on $N$.

To prove  the estimate  (\ref{eq:3.7})  we take  the inner  product  of 
     $Q
\begin{pmatrix}
\frac{d x_\tau}{d\tau}
 \\  
\frac{d \xi_\tau}{d\tau}
\end{pmatrix}$   
with
 $
e^{-2Nt}\begin{pmatrix}
\frac{d x_\tau}{d\tau}
 \\  
\frac{d \xi_\tau}{d\tau}
\end{pmatrix}
$
and  integrate  it in $t$  from  $0$  to $t_0$,   where  $|x_\tau(t_0)|=\max\limits_{0\leq t\leq T}|x_\tau(t)|$.
Note that for any  $\varphi$
\begin{align}																\label{eq:3.8}
&\int\limits_0^{t_0}\frac{d\varphi}{dt}e^{-2Nt}\varphi dt=\frac{1}{2}\int\limits_0^{t_0}e^{-2Nt}\frac{d}{dt}\varphi^2 dt
\\
\nonumber
&=\frac{1}{2}\varphi^2(t_0)e^{-2Nt_0}+N\int\limits_0^{t_0}e^{-2Nt}\varphi^2 dt
\end{align}
Also  we  use in  the proof  of  (\ref{eq:3.7})   that $N$ is large  such  that
\begin{equation}															\label{eq:3.9}
\Big((NI-Q) 
\begin{pmatrix}
\frac{d x_\tau}{d\tau}
 \\  
\frac{d \xi_\tau}{d\tau}
\end{pmatrix}
,
\begin{pmatrix}
\frac{d x_\tau}{d\tau}
 \\  
\frac{d \xi_\tau}{d\tau}
\end{pmatrix}
\Big)>0,
\end{equation}
where  $I$  is the identity operator.
\qed

In addition  to  (\ref{eq:3.7})    we shall estimate  also   $\frac{d^2x_\tau}{d\tau^2}, \frac{d^2\xi_\tau}{d\tau^2}$:      

Differentiating  (\ref{eq:3.3})  in  $\tau$  we get 
\begin{equation}														\label{eq:3.10}
\frac{d}{dt}
\begin{bmatrix} 
\frac{d^2x_\tau}{d\tau^2}
\\
\frac{d^2\xi_\tau}{d\tau^2}
\end{bmatrix}
=
Q
\begin{bmatrix} 
\frac{d^2x_\tau}{d\tau^2}
\\
\frac{d^2\xi_\tau}{d\tau^2}
\end{bmatrix}
+
\frac{dQ}{d\tau}
\begin{bmatrix} 
\frac{dx_\tau}{d\tau}
\\
\frac{d\xi_\tau}{d\tau}
\end{bmatrix}  
+\frac{dF}{d\tau}.
\end{equation}
Therefore
as  in  (\ref{eq:3.6})  we get
\begin{equation}														\label{eq:3.11}
\begin{bmatrix} 
\frac{d^2x_\tau}{d\tau^2}
\\
\frac{d^2\xi_\tau}{d\tau^2}
\end{bmatrix}
=
R(t)
\Big(
\frac{dQ}{d\tau}
\begin{bmatrix} 
\frac{dx_\tau}{d\tau}
\\
\frac{d\xi_\tau}{d\tau}
\end{bmatrix}+\frac{dF}{d\tau}
\Big),
\end{equation}
where  $R(t)$  is the  same  as in  (\ref{eq:3.6}).

Note  that  (cf.  (\ref{eq:3.7}))
\begin{equation}																\label{eq:3.12}
\frac{d Q}{d\tau}=O\Big(\Big|\frac{dx_\tau}{d\tau}\Big|+\Big|\frac{d\xi_\tau}{d\tau}\Big|\Big)
\end{equation}
and
$$
\frac{dF}{d\tau}=
\begin{bmatrix}
\Big((q_2'-q_1')+O\big((q_2'-q_1')^2\big)\Big)\frac{d\xi_\tau}{d\tau}
\\
-\Big(\Big(\frac{\partial q_2'}{\partial x}-\frac{\partial q_1'}{\partial x}\Big)
+O\Big( \frac{\partial q_2'}{\partial x}-\frac{\partial q_1'}{\partial x}\Big)^2\Big)\xi_\tau\frac{d\xi_\tau}{d\tau}
\end{bmatrix}
$$
Since  $\frac{dF}{d\tau}$
can  be  estimated  as in  (\ref{eq:3.7})  we get:
\begin{multline}																\label{eq:3.13}
\max_{0\leq t\leq T}e^{-2Nt}\Big(\Big|\frac{d^2x_\tau}{d\tau^2}\Big| +\Big|\frac{d^2\xi_\tau}{d\tau^2}\Big|              \Big)
\leq C_N \int\limits_0^T e^{-2Nt} \Big(\Big|\frac{dx_\tau}{d\tau}\Big|^2+\Big|\frac{d\xi_\tau}{d\tau}\Big|^2\Big)dt
\\
+C_N\Big(\int\limits_0^T e^{-Nt}|(q_2'-q_1')(x_\tau)|dt\Big)^2+C_N\Big(\int\limits_0^T e^{-Nt}\Big|\frac{\partial}{\partial x}((q_2'-q_1')(x_\tau))\Big|dt\Big)^2.
\end{multline}
Now  we  shall study  the  behavior  in  $\tau$  of
\begin{equation}																\label{eq:3.14}
\frac{dx_0^\tau}{dt}=q_\tau^{00}(x_\tau(t))\xi_0+\sum_{j=1}^nq^{0j}(x_\tau(t))\xi_j(t),\ \ x_0^\tau(0)=0.
\end{equation}
Note that
\begin{equation}																\label{eq:3.15}
q_\tau^{0j}=q_1^{0j}+\tau(q_2^{0j}-q_1^{0j}),\ \ 0\leq j\leq n.
\end{equation}
Therefore
\begin{multline}																	\label{eq:3.16}
\frac{d}{dt}\frac{d}{d\tau}x_0^\tau = \sum_{j=0}^n\big((q_2^{0j}-q_1^{0j})\xi_j  
+O\big(\tau(q_2^{0j}-q_1^{0j})^2\big)\xi_j\big)
+
\sum_{j=0}^n\frac{\partial q_1^{0j}}{\partial  x}\frac{dx^\tau}{d\tau}\xi_j
\\
+ \sum_{j=1}^n q_1^{0j}(x_\tau)\frac{d\xi_j}{d\tau}.
\end{multline}

Thus
\begin{equation}																\label{eq:3.17}
\frac{d}{d\tau}x_0^\tau=\sum_{j=0}^n\int\limits_0^t\Big((q_2^{0j}-q_1^{0j})+O\big(\tau(q_2^{0j}-q_1^{0j})^2\big)\Big)\xi_jdt'+O\Big(\Big|\frac{dx_\tau}{d\tau}\Big|+\Big|\frac{d\xi_\tau}{d\tau}\Big|\Big).
\end{equation}
Denote
\begin{equation}																\label{eq:3.18}
\|q_2^0-q_1^0\|_0=\sum_{j=0}^n\sup_\tau\int\limits_0^T |(q_2^{0j}-q_1^{0j})(x_\tau(t))|dt.
\end{equation}
Then
\begin{equation}																\label{eq:3.19}
\Big|\frac{dx_0^\tau}{d\tau}\Big|\leq C\|q_2^0-q_1^0\|_0 +\max_{0\leq t\leq T}\Big(\Big|\frac{dx_\tau}{d\tau}\Big|+\Big|\frac{d\xi_\tau}{d\tau}\Big|\Big).
\end{equation}

\section{Rigidity  of  the Lorentzian  metric}
\init

Let  $g^{(0)}$  and  $g^{(1)}$  be  two metric  in  $\R\times \Omega$.

Consider  a      $g^{(0)}$  geodesic,  starting at  $y\in \partial\Omega$  when  $t=0$  and  ending  at  $x_T\in\partial\Omega$  when  $t=T$,
and
consider    a     $g^{(1)}$  geodesic,  also starting   at  $y\in\partial\Omega$  when  $t=0$    and ending  at
the same point  $x_T\in\partial\Omega$  when  $t=T$.
The  rigidity property  means that if the lengths  of such two  geodesics  are  equal  and  they  are  close  enough  in some  norm  
then these geodesics  are equal.

Let
\begin{equation}														\label{eq:5.1}
g_\tau=g_0+\tau(g_1-g_0),\ \ \ 0\leq \tau\leq 1.
\end{equation}
Let  $\hat x=\hat x_r(t,\hat y,\hat\eta_\tau)$   be   the  equation  of  the  geodesic  corresponding  to  $g_\tau,  \hat x=(x_0,x)$.
Let  $\hat\eta_0$  be such  that    $ x_0(T, \hat y,\hat\eta_0)=x_T$
and
let  $\hat\eta_1$  be  such  that  $ x_1(T,\hat y,\hat\eta_1)=x_T$.

Consider  the  family  of  $\hat x_\tau(\tau,y,\eta_\tau)$  such  that  $\hat x_\tau$  is  a geodesic  of  $g_\tau$  metric  and
\begin{equation}														\label{eq:5.2}
\hat x_\tau(T,\hat y, \hat\eta_\tau)=x_T\ \ \mbox{for}\ \ 0\leq \tau\leq 1.
\end{equation}
   We have,  differentiating  in $\tau$,
\begin{equation}														\label{eq:5.3}
\frac{d \hat x_\tau}{d\tau}(T,\hat y,\hat\eta_\tau)+\frac{\partial\hat x_\tau}{\partial\eta}\frac{d\hat\eta_\tau}{d\tau}=0.
\end{equation}
Therefore
\begin{equation}														\label{eq:5.4}
\frac{d\hat\eta_\tau}{d\tau}=-\frac{1}{\frac{\partial\hat x_\tau}{\partial\eta}}\frac{d \hat x_\tau}{d\tau},
\end{equation}
assuming that  $\frac{\partial\hat x_\tau}{\partial\eta}\neq 0$.

The  length  of  geodesic  $\hat x_\tau$  is
\begin{multline}															\label{eq:5.5}
R(g_\tau,T,\hat y,\hat\eta_\tau)=
\int\limits_0^T\sqrt{\sum_{j,k=0}^ng_{jk}\frac{dx_j}{dt}\frac{dx_k}{dt}}dt
\\
=\int\limits_0^T\sqrt{\sum_{j,k=0}^ng^{jk}_\tau(x_\tau(t))\xi_j(t)\xi_k(t)}dt
=\sqrt{2H_\tau(y,\hat\eta_\tau)} \, T,
\end{multline}
where
\begin{equation}														\label{eq:5.6}
H_\tau(y,\hat\eta_\tau)=\frac{1}{2}\sum_{j,k=0}^n g_\tau^{jk}(y)\eta_{j\tau}\eta_{k\tau}.
\end{equation}
Note that  
\begin{equation}														\label{eq:5.7}
\frac{dR(g_\tau,T,\hat y,\hat\eta_\tau)}{d\tau}
=(2 H_\tau(y,\hat\eta_\tau))^{-\frac{1}{2}}
\, T\Big( \big((g_1-g_0)\hat\eta_\tau,\hat\eta_\tau\big)
 +2\Big(g_\tau\hat\eta_\tau,\frac{\partial\hat\eta_\tau}{\partial \tau}\Big)\Big)
\end{equation}
and
\begin{equation}														\label{eq:5.8}
\frac{d\hat\eta_\tau}{d\tau}=(-1)\Big(\frac{\partial \hat x_\tau}{\partial\eta}\Big)^{-1}\ \frac{dx_\tau(T,\hat y,\hat\eta_\tau)}{d\tau}.
\end{equation}

We have
\begin{equation}														\label{eq:5.9}
\frac{dR(g_\tau,T,\hat y,\hat\eta_\tau)}{d\tau}\Big|_{\tau=0}=R_1+R_2,
\end{equation}
where
\begin{equation}														\label{eq:5.10}
R_1=(2H(y,\hat\eta_0)^{-\frac{1}{2}}\, T((g_1-g_0)\hat\eta_0,\hat\eta_0),
\end{equation}
\begin{equation}														\label{eq:5.11}
R_2=(2H(y,\hat\eta_0)^{-\frac{1}{2}}\, T\Big(2g_0\hat\eta_\tau,\frac{\partial\hat\eta_\tau}{\partial\tau}\Big)\Big|_{\tau=0}
\end{equation}
Estimating  $\frac{dx_\tau}{d\tau}$  analogously  to  (\ref{eq:3.7})   we get
\begin{equation}														\label{eq:5.12}
\Big|\frac{d\hat\eta_\tau}{d\tau}\Big|\leq C\sup_\tau\int\limits_0^T|(g_1-g_0)(\hat x_\tau(t)\Big|dt
+C\sup_\tau\int\limits_0^T\Big|\frac{\partial}{\partial x}(g_1-g_0)(\hat x_\tau)\Big|dt.
\end{equation}

Denote
\begin{equation}														\label{eq:5.13}
\||(g_1-g_0)\||=\sup_\tau\int\limits_0^T|(g_1-g_0)(\hat x_\tau)|dt
+\sup_\tau\int\limits_0^T\Big|\Big(\frac{\partial g_1}{\partial x}-\frac{\partial g_0}{\partial x}\Big)(\hat x_\tau)\Big|dt
\end{equation}
and denote 
\begin{equation}														\label{eq:5.14}
\|g_1-g_0\|=\sup_y|g_1(y)-g_0(y)|.
\end{equation}
Choose  $w=\frac{g_1-g_0}{\|g_1-g_0\|}$     such  that  
\begin{equation}														\label{eq:5.15}
R_1=l_0\|g_1-g_0\|, \ \ \mbox{where}\ \ \
l_0=(2H(y,\hat\eta_0))^{-\frac{1}{2}}\, T(w\hat\eta_0,\hat\eta_0)>0.
\end{equation}
Analogously    we have
\begin{equation}														\label{eq:5.16}
R_2=(2H(y,\hat\eta_0))^{-\frac{1}{2}}\, T\Big(2g_0\hat\eta_1,\frac{\partial\hat\eta_1}{\partial\tau}\Big|_{\tau=0}\Big)=
\alpha(g_1-g_0),
\end{equation}
where $\alpha(g_1-g_0)$  is  a  linear  functional.

Denote  
\begin{equation}														\label{eq:5.17}
w_1=\frac{g_1-g_0}{|||g_1-g_0|||},
\end{equation}
where   the norm  $|||g_1-g_0|||$  is the same  as  in  (\ref{eq:5.13}).

We  choose $w_1$  without  changing  $w$  such that
\begin{equation}														\label{eq:5.18}
\alpha(g_1-g_0)=l_1|||g_1-g_0|||,
\end{equation}
where  $l_1=\alpha(w_1)>0$.

Therefore   we get
\begin{equation}														\label{eq:5.19}
R_1+R_2=l_0\|g_1-g_0\|+l_1|||g_1-g_0|||>0.
\end{equation}
Let
\begin{equation}														\label{eq:5:20}
\Delta=R(g_1,T,\hat y,\hat\eta_1)-   R(g_0,T,\hat y,\hat \eta_0).
\end{equation}  
We have
\begin{multline}														\label{eq:5.21}
|R(g_2,T,y,\eta_1)-R(g_0,T,y,\eta_0)-R_1-R_2|
\\
\leq C_1\|g_1-g_0\|^2+C_2|||g_1-g_0|||^2.
\end{multline}
Therefore
\begin{equation}												\label{eq:5.22}
R_1+R_2\leq |\Delta|+C_1\|g_1-g_0\|^2+C_2|||g_1-g_0|||^2.
\end{equation}

If  $\|g_1-g_0\|$  and  $|||g_1-g_0|||$  are small  enough,  more precisely,   if
\begin{equation}														
\nonumber
\|g_1-g_0\|<\frac{1}{2l_0C_1}
\ \ \ \ \mbox{and}\ \ \ \ \
|||g_1-g_0|||<\frac{1}{2l_1 C_2},
\end{equation}
then
\begin{equation}														
\nonumber
\frac{1}{2}l_0\|g_1-g_0\|\leq|\Delta|
\ \ \ \ \mbox{and}\ \ \ \ \
\frac{1}{2}l_1|||g_1-g_0|||\leq|\Delta|.
\end{equation}
From   the last  inequalities  we get  that  $|\Delta|=0$  implies  that  $g_1-g_0=0$.

Thus  we proved  the following  theorem:
\begin{theorem}														\label{theo:5.1}
Let  $R(g_1,T,y,\hat\eta_1)=R(g_0,T,y,\hat\eta_0)$,   where  $\hat\eta_1$  and  $\hat\eta_0$  are  such  that 
$x_0(T,y,\hat\eta_0)=x_T,
x_1(T,y,\hat\eta_1)=x_T$,  i.e.  the geodesics  in both  metrics  have the same length.
If metrics $g_0$  and  $g_1$  are close  in norms  (\ref{eq:5.13})  and  (\ref{eq:5.14})  then  $g_0=g_1$.
\end{theorem}

\section{Rigidity in the case  of  null-geodesics}
\init
The length  of  space-time null-geodesics  
$x_0=x_0(t,y,\eta),\  x=x(t,y,\eta),\linebreak 0\leq t\leq T,$   where   $x(T,y,\eta)\in\partial\Omega$,  is equal  to
\begin{equation}																\label{eq:4.1}
L(q,T,y,\eta)=\int\limits_0^T\sqrt{\Big(\frac{dx_0}{dt}\Big)^2+\Big(\frac{dx(q')}{dt}\Big)^2}dt,
\end{equation}
where  $x(q')(t)$  is the solution  of Hamiltonian system
$$
\frac{dx(q',t,y,\eta)}{dt}=q'(x(t))\xi(t),\ \ \frac{d\xi(q',t,y,\eta)}{dt}=-\frac{1}{2}\Big(\frac{\partial q'(x)}{\partial x}\xi(t)\Big)\cdot\xi(t),
$$
$x\big|_{t=0}=y,\ \xi\big|_{t=0}=\eta,\ 0\leq t\leq  T,\ \ \big|\frac{dx}{dt}\big|=\sqrt{\sum_{k=1}^n\big(\frac{dx_k}{dt}\big)^2},\\ \frac{dx_0}{dt}=
\sum_{j=1}^n q^{0j}(x(t))\xi_j(t)+q^{00}(x(t))\eta_0,\ x_0(0)=0$.

Let  $q=q_1+\tau(q_2-q_1),  0\leq \tau\leq 1$.  We have
\begin{multline}																\label{eq:4.2}
\frac{\partial L(q,T,y,\eta)}{d\tau}=\int\limits_0^T\Big( \Big(\frac{dx_0}{dt}\Big)^2+\Big|\frac{dx}{dt}\Big|^2 \Big)^{-\frac{1}{2}}\Big(\frac{dx}{dt},\frac{d}{d\tau}\frac{dx}{dt}\Big)dt
\\
+\int\limits_0^T\Big(\Big(\frac{dx_0}{dt}\Big)^2+\Big|\frac{dx}{dt}\Big|^2\Big)^{-\frac{1}{2}}\Big(\frac{dx_0}{dt},\frac{d}{d\tau}\frac{dx_0}{dt}\Big)dt.
\end{multline}
Note that
$$
\frac{d}{d\tau}\frac{dx}{dt}=\frac{d}{d\tau}(q'\xi)=   (q_2'-q_1')\xi
+ \tau \frac{d}{d\tau}\Big( (q_2'-q_1') \xi      \Big)   
+\frac{\partial (q_2'-q_1')(x(t))}{\partial x}\ \frac{dx}{d\tau}\xi  +q_1'(x)\frac{d\xi}{d\tau}.
$$
Therefore
\begin{multline}																\label{eq:4.3}
\frac{\partial L}{\partial \tau}\Big|_{\tau=0}
=\int\limits_0^T
\Big( \Big(\frac{dx_0}{dt}\Big)^2+\Big|\frac{dx}{dt}\Big|^2 \Big)^{-\frac{1}{2}}
\Big(q_1'(x(t))\xi(t),(q_2'-q_1')\xi
+\frac{\partial q_1'}{\partial x}\xi \frac{dx}{d\tau}+q_1'(x)\frac{d\xi}{d\tau}\Big)
\\
+\Big(\frac{dx_0}{dt},\frac{d}{d\tau}\frac{dx_0}{dt}\Big|_{t=0}\Big)\Big)
dt.
 \end{multline}
Thus  by the Taylor's  formula
\begin{equation}																\label{eq:4.4}
L(q_2,T,y,\eta)-L(q_1,T,y,\eta)=\tau\frac{\partial L(q,T,y,\eta)}{\partial \tau}\Big|_{\tau=0}+G_2,
\end{equation}
where
\begin{equation}																\label{eq:4.5}
G_2=\frac{1}{2}\frac{\partial^2}{\partial\tau^2}L(q_{1}+	\theta(q_2-q_1),T,y,\eta)(q_2-q_1)^2,\ 0<\theta <1.
\end{equation}
Note  that  
\begin{equation}																\label{eq:4.6}
\tau\frac{\partial L(q_1,y,T,\eta)}{\partial\tau}\big|_{\tau=0}=l(\tau(q_2-q_1))																			
\end{equation}
is  the  linear  part  of  $L(q_2)-L(q_1)$.

Let
\begin{multline}																\label{eq:4.7}
\|q_2-q_1\|=\sup_\tau\int\limits_0^T e^{-2Nt}|(q_2'-q_1')(x_\tau(t)|dt
\\
+\sup_\tau\int\limits_0^T  e^{-2Nt}\Big|\frac{\partial}{\partial x}(q_2'-q_1')(x_\tau(t))\Big|dt
+\|q_2^0-q_1^0\|_0,
\end{multline}
where    $\|q_2^0-q_1^0\|_0$  is the same  as   in (\ref{eq:3.18}).
Let  $\omega(q_2-q_1)= \frac{q_2-q_1}{\|q_2-q_1\|}$.
Since  $l(q_2-q_1)$  is  nonzero linear  functional  bounded  in the norm  (\ref{eq:4.7}) 
there exists  $\omega_1$  such that  $|l(\omega_1)|=l_0>0$.    Therefore
\begin{equation}																\label{eq:4.8}
|l(q_2-q_1)|=|l(\omega_1)|\, \|q_2-q_1\|=l_0\|q_2-q_1\|.
\end{equation}
Now  estimate $G_2$. 
Denote 
$$
\lambda= \Big(\frac{dx}{dt}\Big)^2 +\Big(\frac{dx_0}{dt}\Big)^2.
$$
 Differentiating  $L(q_1+\tau(q_2-q_1))$  twice  in $\tau$  we  get
\begin{align}																	\label{eq:4.9}
G_2&= \frac{d^2}{d\tau^2}    \int\limits_0^T \sqrt\lambda dt
\\
\nonumber
&=\int\limits_0^T\Big[\Big(-\frac{1}{4}\Big)\lambda^{-{3}/{2}}\Big(\frac{d\lambda}{d\tau}\Big)^2
+\frac{1}{2}\lambda^{-{1}/{2}}
\frac{d^2\lambda}{d\tau^2}\Big]dt
\end{align}

Estimating the right hand sides  in (\ref{eq:4.9})  as in   (\ref{eq:3.7}),  (\ref{eq:3.13}),  (\ref{eq:3.19})   we get
\begin{equation}																\label{eq:4.10}
|G_2|\leq  C\int\limits_0^T\Big[ \Big(\Big|\frac{d^2 x}{d\tau^2}\Big|+\Big|\frac{d^2\xi}{d\tau^2}\Big|+\Big|\frac{d^2x_0}{d\tau^2}\Big|\Big)
+C\Big(\Big|\frac{dx}{d\tau}\Big|^2+\Big|\frac{d\xi}{d\tau}\Big|^2+\Big(\frac{dx_0}{d\tau}\Big)^2\Big) \Big]
dt.
\end{equation}
 Using  (\ref{eq:3.13})  and  (\ref{eq:3.19})  we obtain
\begin{equation}																\label{eq:4.11}
|G_2|\leq  C_N \|q_2-q_1\|^2.
\end{equation}
Since
\begin{equation}																\label{eq:4.12}
L(q_2,T,y,\eta)-L(q_1,T,y,\eta)
=l(q_2-q_1)+G_2,
\end{equation}
we  have, assuming that $\omega=\omega_1$  and  using  (\ref{eq:4.10})  and  (\ref{eq:4.11}):
\begin{equation}																\label{eq:4.13}
l_0\|q_2-q_1\|\leq |L(q_2,T,y,\eta)-L(q_1,T,y,\eta)|+C_N\|q_2-q_1\|^2.
\end{equation}	
Therefore
\begin{equation}																\label{eq:4.14}
l_0\|q_2-q_1\|\Big(1-\frac{C_N}{l_0}\|q_2-q_1\|\Big)\leq |L(q_2,T,y,\eta)-L(q_1,T,y,\eta)|.
\end{equation}
Assuming  that  $\|q_2-q_1\|<\frac{l_0}{2C_N}$  we obtain
\begin{equation}																\label{eq:4.15}
2l_0\|q_2-q_1\|\leq  |L(q_2,T,y,\eta)-L(q_1,T,y,\eta)|.
\end{equation}
Then   $L(q_2,T,y,\eta)=L(q_1,T,y,\eta)$  implies  that  $\|q_2-q_1\|=0$. 

Let $x'_0(t)$  be  the  null-geodesic  in  $q_1'$   metric  starting  at  $y_0$  when  $t=0$,  
such that it  is
 perpendicular 
to $\partial\Omega$
 at $y_0$.
 We have:  
$$\int\limits_0^T\big| \big(q'_2(x'_0(t))-q'_1(x'_0(t)\big)\big|dt=0,$$
  which implies that  $q'_2(x'_0(t))-q'_1(x'_0(t)=0$ a.e. on $[0,T]$.
  We have also that
$q_1^0(x'_0(t))=q_2^0(x'_0(t))$. 
Thus  $q_1(x'_0(t))=q_2(x'_0(t))$.
Since this is true for all $x'_0(t)$  
we get that
$q_1(x)=q_2(x)$  on $\Omega$.

  Let  $T=T(y,\eta)$  be  such  that  $x_T=x(T,y,\eta)\in  \partial\Omega$,  
i.e.   $x(t,y,\eta)$  leaves  $\R\times\Omega$  when  $t>T(y,\eta)$.  We shall call such  $T(y,\eta)$  maximal.

The following theorem  holds:
\begin{theorem}															\label{theo:4.1}
Let  $\sum_{j,k=0}^n q_{pjk}(x)dx_j dx_k$  be  two metrics,  $p=1,2,$  
and let  $(x_0(q_p),x(q_p)))$  be  the space-time  null-geodesics with the same  initial  conditions  $x_0=0,  x=y$  and  $\xi=\eta$.
 Let $T'(y,\eta)$  be  maximal  in  $q_1$-metric   for all  $y$.
Then  if  $L(q_2,T'(y,\eta),y,\eta)=L(q_1,T'(y,\eta),y,\eta)$  for  all $y$
and  if   $q_2$  and  $q_1$   are sufficiently   close (as in (\ref{eq:4.15})), 
 then  $q_2=q_1$.

\end{theorem}

{}

\end{document}